\documentclass[a4paper,12pt]{amsart}
\usepackage{amssymb,amscd,amsmath,a4wide,graphicx,stmaryrd,fullpage,setspace,microtype,xr,textcomp,enumitem,wasysym,centernot}
\usepackage{hyperref}
\usepackage[all]{xy}
\usepackage[T1]{fontenc}
\usepackage{lmodern}
\usepackage[numbers,sort&compress]{natbib}

\title{Comparing cubical and globular directed paths}

\author[P. Gaucher]{Philippe Gaucher}

\address{Universit\'e Paris Cit\'e, CNRS, IRIF, F-75013, Paris, France}

\urladdr{http://www.irif.fr/{\~{}}gaucher} 

\makeatletter
\@namedef{subjclassname@2020}{\textup{2020} Mathematics Subject Classification}
\makeatother
\subjclass[2020]{55U35,68Q85}

\keywords{directed path, precubical set, directed homotopy, Reedy category, combinatorial model category, accessible model category}

\swapnumbers

\newcommand{\C}{\mathcal{C}}
\newcommand{\D}{\mathcal{D}}

\newcommand{\de}{\partial}
\newcommand{\p}{\times}

\renewcommand{\P}{\mathbb{P}}

\newtheorem*{thmN}{Theorem}

\newtheorem{thm}{Theorem}[section]
\newtheorem{prop}[thm]{Proposition}

\newtheorem{exa}[thm]{Example}
\newtheorem{cor}[thm]{Corollary}

\newtheorem{defn}[thm]{Definition}
\newtheorem{nota}[thm]{Notation}
\newtheorem{rem}[thm]{Remark}
\newcommand{\bd}{\begin{defn}}
	\newcommand{\ed}{\end{defn}}
\newcommand{\bp}{\begin{prop}}
	\newcommand{\ep}{\end{prop}}
\newcommand{\bth}{\begin{thm}}
	\renewcommand{\eth}{\end{thm}}
\newcommand{\bpf}{\begin{proof}}
	\newcommand{\epf}{\end{proof}}
\newcommand{\bc}{\begin{cor}}
	\newcommand{\ec}{\end{cor}}

\newcommand{\fL}[1]{\ar@{->}[ll]_-{#1}}
\newcommand{\fR}[1]{\ar@{->}[rr]^-{#1}}
\newcommand{\fRr}[1]{\ar@{->}[rrr]^-{#1}}
\newcommand{\fD}[1]{\ar@{->}[dd]_-{#1}}
\newcommand{\fU}[1]{\ar@{->}[uu]^-{#1}}
\newcommand{\f}[2]{\ar@{->}[#1]|{#2}}
\newcommand{\ff}[2]{\ar@2{->}[#1]|{#2}}
\newcommand{\frr}[1]{\ar@{->}[rrrr]^-{#1}}

\newcommand{\fl}[1]{\ar@{->}[l]_-{#1}}
\newcommand{\fr}[1]{\ar@{->}[r]^-{#1}}
\newcommand{\fd}[1]{\ar@{->}[d]_-{#1}}
\newcommand{\fu}[1]{\ar@{->}[u]^-{#1}}

\renewcommand{\top}{{\mathbf{Top}}}

\newcommand{\iso}{\cong}

\renewcommand{\leq}{\leqslant}
\renewcommand{\geq}{\geqslant}

\DeclareMathOperator{\carrier}{Carrier}

\def\cartesien{%
	\ar@{-}[]+R+<6pt,-2pt>;[]+RD+<6pt,-6pt>%
	\ar@{-}[]+D+<2pt,-6pt>;[]+RD+<6pt,-6pt>%
}
\def\cocartesien{%
	\ar@{-}[]+L+<-6pt,+2pt>;[]+LU+<-6pt,+6pt>%
	\ar@{-}[]+U+<-2pt,+6pt>;[]+LU+<-6pt,+6pt>%
}
\def\hocartesien{%
	\ar@{-}[]+R+<6pt,-2pt>;[]+RD+<6pt,-6pt>_{h}%
	\ar@{-}[]+D+<2pt,-6pt>;[]+RD+<6pt,-6pt>%
}
\def\hococartesien{%
	\ar@{-}[]+L+<-6pt,+2pt>;[]+LU+<-6pt,+6pt>_{h}%
	\ar@{-}[]+U+<-2pt,+6pt>;[]+LU+<-6pt,+6pt>%
}

\newcommand{\brm}[1]{{\rm{\mathbf{#1}}}}

\newcommand{\dtop}{{\brm{Flow}}}

\newcommand{\set}{{\brm{Set}}}

\newcommand{\ttop}{{\brm{TOP}}}

\newcommand{\glob}{{\mathrm{Glob}}}

\DeclareMathOperator{\id}{Id}

\DeclareMathOperator{\Ch}{Ch}

\newcommand{\mins}[1]{\texorpdfstring{$#1$}{Lg}}

\newcommand{\liminj}{\varinjlim}

\renewcommand{\P}{\mathbb{P}}

\makeatletter
\def\varholim@#1#2{%
	\vtop{\m@th\ialign{##\cr
			\hfil$#1\operator@font holim$\hfil\cr
			\noalign{\nointerlineskip\kern1.5\ex@}#2\cr
			\noalign{\nointerlineskip\kern-\ex@}\cr}}%
}
\def\holimproj{%
	\mathop{\mathpalette\varholim@{\leftarrowfill@\textstyle}}\nmlimits@
}
\def\holiminj{%
	\mathop{\mathpalette\varholim@{\rightarrowfill@\textstyle}}\nmlimits@
}
\makeatother

\newcommand{\ddownarrow}{{\downarrow}}

\newcommand{\adj}[5]{\xymatrix@C=#5em{{#1}\ar@/^0.8em/[r]^-{#2} \ar@{}[r]|-{\perp} & \ar@/^0.8em/[l]^-{#3} {#4}}}

\DeclareMathOperator{\seq}{Seq}

\allowdisplaybreaks

\begin{document}

\begin{abstract}
	A flow is a directed space structure on a homotopy type. It is already known that the underlying homotopy type of the realization of a precubical set as a flow is homotopy equivalent to the realization of the precubical set as a topological space. This realization depends on the non-canonical choice of a q-cofibrant replacement. We construct a new realization functor from precubical sets to flows which is homotopy equivalent to the previous one and which does not depend on the choice of any cofibrant replacement functor. The main tool is the notion of natural $d$-path introduced by Raussen. The flow we obtain for a given precubical set is not anymore q-cofibrant but is still m-cofibrant. As an application, we prove that the space of execution paths of the realization of a precubical set as a flow is homotopy equivalent to the space of nonconstant $d$-paths between vertices in the geometric realization of the precubical set. 
\end{abstract}
	
\maketitle
\tableofcontents

\section{Introduction}

\subsection*{Presentation}

Precubical sets are a prominent geometric model for concurrency theory \cite{DAT_book}. The $n$-cube represents the concurrent execution of $n$ actions. The space of $d$-paths in the geometric realization of a precubical set is studied in many papers, such as a series of papers \cite{MR2521708,Raussen2010,Raussen2012,MR3722069,MR4070250} by Raussen and Ziemia\'{n}ski (the list of references is not exhaustive). Precubical sets can also be realized as flows in the sense of \cite{model3}. The realization functor of a precubical set as a flow is first introduced in \cite[Definition~7.2]{ccsprecub}.  

The two approaches (let us call them the cubical one of Raussen and Ziemia\'{n}ski and the globular one of the author) do not coincide up to homeomorphism. In the cubical approach, the $d$-paths of the topological $n$-cube $[0,1]^n$ are the continuous paths which are nondecreasing with respect to each axis of coordinates. Raussen and Ziemia\'{n}ski study also several variants (tame, strict, natural etc...) which give rise to homotopy equivalent spaces of $d$-paths between two fixed vertices in the geometric realization of the precubical set. None of these definitions give rise to spaces of $d$-paths from the initial to the final states of the $n$-cube $[0,1]^n$ which are homeomorphic to the spaces of execution paths from the initial to the final states of the $n$-cube $[0,1]^n$ viewed as a flow. In the latter case, the space of execution paths from the initial to the final states of the $n$-cube is the $(n-1)$-dimensional disk $\mathbf{D}^{n-1}$ (see Theorem~\ref{rea-func}). It means that the latter space depends on a non-canonical choice of an achronal slice in the middle  of the topological $n$-cube and on a non-canonical choice of a homeomorphism between this achronal slice and $\mathbf{D}^{n-1}$.

The underlying homotopy type of a flow is the homotopy type obtained after removing the temporal information contained in a flow. It is defined in \cite[Section~6]{4eme} and a more conceptual construction is provided in \cite[Proposition~8.16]{Moore2} using Moore flows. It is already known in full generality that the underlying homotopy type of the realization of a precubical set as a flow is isomorphic to the homotopy type of the realization of a precubical set as a topological space \cite[Theorem~6.2.1]{realization}. One of the purposes of this paper is to prove the directed version of this result.

At first, using the notion of natural $d$-path introduced by Raussen in \cite[Definition~2.14]{MR2521708}, we improve the realization functor from precubical sets to flows $|\!-\!|_q:\square^{op}\set\to \dtop$ introduced in \cite[Definition~7.2]{ccsprecub} as follows.

\begin{thmN} (Theorem~\ref{natural-equal})
	There exist a colimit preserving functor \[|\!-\!|_{nat}:\square^{op}\set\longrightarrow \dtop\] which does not depend on any cofibrant replacement and a natural transformation $\mu:|\!-\!|_q\Rightarrow |\!-\!|_{nat}$ such that for all precubical sets $K$, the natural map $\mu_K:|K|_q \to |K|_{nat}$ induces a bijection on states and a homotopy equivalence $\P_{\alpha,\beta}|K|_q \simeq \P_{\alpha,\beta}|K|_{nat}$ for all $\alpha,\beta\in K_0$.
\end{thmN}

Theorem~\ref{natural-equal} implies that $\P_{\alpha,\beta}|K|_{nat}$ is m-cofibrant, $\P_{\alpha,\beta}|K|_q$ being q-cofibrant. The interest of the natural realization functor is that it does not depend anymore on the arbitrary choice on any cofibrant replacement functor for the category of flows. Surprisingly, it does not even depend on an m-cofibrant replacement or on an h-cofibrant replacement of the category of flows. The geometric properties of the natural $d$-paths enable us in Section~\ref{cube_sec} to give another description of the natural realization functor using Ziemia{\'{n}}ski's notion of cube chain. As an application of Theorem~\ref{natural-equal} and of Section~\ref{cube_sec}, we prove the following theorems:

\begin{thmN} (Theorem~\ref{main2} and Theorem~\ref{final})
 Let $K$ be a precubical set. Let $\alpha,\beta$ be two vertices of $K$. The space of execution paths $\P_{\alpha,\beta}|K|_{nat}$ is homotopy equivalent to the set of nonconstant tame natural $d$-paths from $\alpha$ to $\beta$ equipped with the $\Delta$-kelleyfication of the relative topology induced by the compact-open topology in the geometric realization of $K$. When moreover $K$ is spatial (e.g. proper), the homotopy equivalence is a homeomorphism. 
\end{thmN}

\begin{thmN} (Corollary~\ref{app_general})
	Let $K$ be a precubical set. Let $\alpha,\beta$ be two vertices of $K$. Then the space of execution paths $\P_{\alpha,\beta}|K|_q$ is homotopy equivalent to the space of nonconstant $d$-paths from $\alpha$ to $\beta$ in the geometric realization of $K$ equipped with the $\Delta$-kelleyfication of the compact-open topology. 
\end{thmN}

Corollary~\ref{app_general} is not surprising. However, until a proof was known, it was not sure that the statement was true for all precubical sets in full generality, and not only e.g. for non-positively curved precubical sets in the sense of \cite{zbMATH07226006}, notion which brings together the properties satisfied by the precubical sets coming from a lot of real concurrent systems by \cite[Proposition~1.29]{zbMATH07226006}. 

\subsection*{Outline of the paper}

Section~\ref{three_section} is a reminder about the three model structures of flows: Quillen (q), Hurewicz (h) and mixed (m) introduced in \cite{QHMmodel}. It contains, as a new and easy remark, the proof that these three model structures on flows are simplicial. Section~\ref{rea_sec} recalls some basic facts about cocubical objects, gives the definition of an r-realization functor with $r\in \{q,m,h\}$ in Definition~\ref{def-rea-flow}, adapts in Theorem~\ref{homotopy-natural} some tools coming from \cite{realization}, and finally gives the example of the q-realization functor expounded in \cite{ccsprecub}. Section~\ref{natural_sec} recalls the notion of tame natural $d$-path of a precubical set and proves some basic facts about their topology, in relation with the $\Delta$-generated spaces which are the setting of this work. Section~\ref{nat_rea_sec} expounds the construction of the natural realization functor in Definition~\ref{natural-rea}. It does not depend on any cofibrant replacement functor. It is a new realization functor which is proved to be equivalent in some sense to the one of \cite{ccsprecub} in Theorem~\ref{natural-equal}. This section also proves that this new realization functor is an m-realization functor. Section~\ref{cube_sec} gives an equivalent definition in Theorem~\ref{iso-cube} of the natural realization functor in terms of cube chains in the sense of Ziemia{\'{n}}ski. In Section~\ref{spatial_precubical_set}, the tame concrete realization of a precubical set as a flow and the notion of spatial precubical set are introduced in Definition~\ref{rea_conc} and in Definition~\ref{spatial} respectively. Theorem~\ref{main2} proves that the natural realization and the tame concrete realization coincide in the spatial case. Then the latter theorem is generalized to the non-spatial case in Theorem~\ref{final}. Finally it is described the connexion with the q-realization functor in Corollary~\ref{app_general}. Appendix~\ref{spatial_app} is devoted to the proof that the class of spatial precubical sets is a small orthogonality class.

\subsection*{Prerequisites and notations}

We refer to \cite{TheBook} for locally presentable categories, to \cite{MR2506258} for combinatorial model categories.  We refer to \cite{MR99h:55031} and to \cite{ref_model2} for more general model categories, and to \cite{zbMATH06722019,HKRS17,GKR18} for accessible model categories. The main tools used in this paper are the $\{q,m,h\}$-model structures of flows \cite{QHMmodel}, the homotopical results of \cite{realization} about the realization functors of precubical sets as flows, and some topological results due to Ziemia{\'{n}}ski about natural $d$-paths and the technique of cube chains coming from \cite{MR4070250}. The initial object of a category is denoted by $\varnothing$. The terminal object of a category is denoted by $\mathbf{1}$. The set of maps from $X$ to $Y$ in a category $\C$ is denoted by $\C(X,Y)$. $\id_X$ denotes the identity map of $X$. The category $\top$ denotes the category of \textit{$\Delta$-generated spaces} or of \textit{$\Delta$-Hausdorff $\Delta$-generated spaces} (cf. \cite[Section~2 and Appendix~B]{leftproperflow}). The inclusion functor from the full subcategory of $\Delta$-generated spaces to the category of general topological spaces together with the continuous maps has a right adjoint called the $\Delta$-kelleyfication functor. The latter functor does not change the underlying set. The category $\top$ is locally presentable and cartesian closed. The internal hom $\ttop(X,Y)$ is given by taking the $\Delta$-kelleyfication of the compact-open topology on the set $\top(X,Y)$. The category $\top$ is equipped with its q-model structure. The m-model structure \cite{mixed-cole} and the h-model structure \cite{Barthel-Riel} of $\top$ are also used in various places of the paper. The q-model structure of $\top$ is combinatorial. The m-model structure and the h-model structure of $\top$ are accessible~\footnote{It is unlikely that they are combinatorial but no proof is known. The proof of \cite[Remark~4.7]{Raptis2} that the h-model structure is not cofibrantly generated works only for the category of \textit{general} topological spaces.}. The three model structures are monoidal and simplicial. \textit{Compact means quasicompact Hausdorff (French convention)}. 

\subsection*{Warning}

All $d$-paths in a geometric realization of a precubical set considered in this paper are tame and nonconstant: see Remark~\ref{warning}. The adjective tame is added on purpose everywhere. The adjective nonconstant is often omitted (but always understood), except in Corollary~\ref{app_general} to avoid any confusion.

\subsection*{Acknowledgment}

I thank Krzysztof Ziemia{\'{n}}ski for helpful discussions about \cite{MR4070250}. I am grateful to the anonymous referee for pointing out the flaw in the proof of the main theorem, which also led me to many improvements in the exposition of the results.

\section{Three simplicial model structures of flows}
\label{three_section}

\bd \cite[Definition~4.11]{model3} \label{def-flow}
A {\rm flow} is a small semicategory enriched over the closed monoidal category $(\top,\p)$. The corresponding category is denoted by $\dtop$. 
\ed

A \textit{flow} $X$ consists of a topological space $\P X$ of \textit{execution paths}, a discrete space $X^0$ of \textit{states}, two continuous maps $s$ and $t$ from $\P X$ to $X^0$ called the source and target map respectively, and a continuous and associative map $*:\{(x,y)\in \P X\p \P X; t(x)=s(y)\}\longrightarrow \P X$ such that $s(x*y)=s(x)$ and $t(x*y)=t(y)$. Let $\P_{\alpha,\beta}X = \{x\in \P X\mid s(x)=\alpha \hbox{ and } t(x)=\beta\}$: it is the space of execution paths from $\alpha$ to $\beta$, $\alpha$ is called the initial state and $\beta$ is called the final state. Note that the composition is denoted by $x*y$, not by $y\circ x$. The category $\dtop$ is locally presentable by \cite[Theorem~6.11]{Moore1}.

Every set can be viewed as a flow with an empty space of execution paths. Every poset can be viewed as a flow with one execution path from $\alpha$ to $\beta$ if and only if $\alpha<\beta$. The obvious functor $\set \subset \dtop$ from the category of sets to that of flows is limit-preserving and colimit-preserving. The following example of flows is important for the sequel: 

\begin{exa}
	For a topological space $Z$, let $\glob(Z)$ be the flow defined by 
	\[
	\glob(Z)^0=\{0,1\}, \ 
	\P \glob(Z)= \P_{0,1} \glob(Z)=Z,\ 
	s=0,\  t=1.
	\]
	This flow has no composition law.
\end{exa}

\begin{nota}
	Let $n\geq 1$. Denote by $\mathbf{D}^n = \{b\in \mathbb{R}^n, |b| \leq 1\}$ the $n$-dimensional disk, and by $\mathbf{S}^{n-1} = \{b\in \mathbb{R}^n, |b| = 1\}$ the $(n-1)$-dimensional sphere. By convention, let $\mathbf{D}^{0}=\{0\}$ and $\mathbf{S}^{-1}=\varnothing$.
\end{nota}

We need to recall:

\bth \label{three} Let $r\in \{q,m,h\}$. Then there exists a unique model structure on $\dtop$ such that: 
\begin{itemize}
	\item A map of flows $f:X\to Y$ is a weak equivalence if and only if $f^0:X^0\to Y^0$ is a bijection and for all $(\alpha,\beta)\in X^0\p X^0$, the continuous map $\P_{\alpha,\beta}X\to \P_{f(\alpha),f(\beta)}Y$ is a weak equivalence of the r-model structure of $\top$.
	\item A map of flows $f:X\to Y$ is a fibration if and only if for all $(\alpha,\beta)\in X^0\p X^0$, the continuous map $\P_{\alpha,\beta}X\to \P_{f(\alpha),f(\beta)}Y$ is a fibration of the r-model structure of $\top$.
\end{itemize}
This model structure is accessible and all objects are fibrant. Moreover, this model structure is simplicial. It is called the r-model structure of $\dtop$.
\eth

\bpf It is \cite[Theorem~7.4]{QHMmodel} except the last sentence. The fact that $\dtop$ is enriched, tensored and cotensored over simplicial sets is proved in \cite[Section~3.3]{realization}. The q-model structure of flows is simplicial by \cite[Theorem~3.3.15]{realization}. It remains to prove the compatibility with the m-model structure and the h-model structure. It suffices to prove (see the very end of the proof of \cite[Proposition~3.3.14]{realization}) that the lift
$k'$ of the commutative square
\[
\xymatrix{ (\mathbf{D}^n\p \{0\}\p |\Delta[1]|)\cup (\mathbf{D}^n\p[0,1]\p
	\{-1,1\}) \fr{}\fd{} &
	\P_{\alpha,\beta}X\fd{}\\
	\mathbf{D}^n\p[0,1]\p |\Delta[1]| \fr{} \ar@{-->}[ru]_-{k'}& \P_{\alpha,\beta}Y}
\] 
exists if the map $\P_{\alpha,\beta}X\to \P_{\alpha,\beta}Y$ is an m-fibration of spaces or an h-fibration of spaces. Since every m-fibration of spaces and every h-fibration of spaces is a q-fibration of spaces, the proof is complete.
\epf

By \cite[Theorem~7.7]{QHMmodel}, the m-model structure is the mixing of the q-model structure and the h-model structure in the sense of \cite[Theorem~2.1]{mixed-cole}. The q-model structure is not only accessible, but also combinatorial. A set of generating cofibrations is the set of maps $\{\glob(\mathbf{S}^{n-1})\subset \glob(\mathbf{D}^{n}) \mid n\geq 0\} \cup \{C:\varnothing \to \{0\},R:\{0,1\} \to \{0\}\}$ by e.g. \cite[Theorem~7.6]{QHMmodel}. Every q-cofibration of flows is an m-cofibration and every m-cofibration of flows is an h-cofibration by \cite[Proposition~3.6]{mixed-cole}. 

There exists a flow which is not cofibrant in any of the three model structures by \cite[Proposition~7.9]{QHMmodel}. This behaviour differs from the behaviour of the h-model structure of topological spaces for which all spaces are h-cofibrant (and h-fibrant). The reason is that the h-model structure of flows does not coincide with the Hurewicz model structure given by \cite[Corollary~5.23]{Barthel-Riel}. This one exists as well because $\dtop$ satisfies the monomorphism hypothesis, being locally presentable, and because $\dtop$ is topologically bicomplete (the proof is similar to the proof that it is simplicial as given in \cite[Section~3.3]{realization}) since a $\Delta$-generated space is homeomorphic to the disjoint sum of its path-connected components by \cite[Proposition~2.8]{mdtop}. This Hurewicz model structure is not used in this paper.

\section{Realization functors from precubical sets to flows}
\label{rea_sec}

\begin{nota}
	Let $[0] = \{()\}$ and $[n] = \{0,1\}^n$ for $n \geq 1$. By convention, one has $\{0,1\}^0=[0]=\{()\}$. The set $[n]$ is equipped with the product ordering $\{0<1\}^n$. Let $0_n=(0,\dots,0) \in \{0,1\}^n$ and $1_n=(1,\dots,1) \in \{0,1\}^n$
\end{nota}

Let $\delta_i^\alpha : [n-1] \rightarrow [n]$ be the coface map defined for $1\leq i\leq n$ and $\alpha \in \{0,1\}$ by $\delta_i^\alpha(\epsilon_1, \dots, \epsilon_{n-1}) = (\epsilon_1,\dots, \epsilon_{i-1}, \alpha, \epsilon_i, \dots, \epsilon_{n-1})$. The small category $\square$ is by definition the subcategory of the category of sets with the set of objects $\{[n],n\geq 0\}$ and generated by the morphisms $\delta_i^\alpha$. They satisfy the cocubical relations $\delta_j^\beta \delta_i^\alpha = \delta_i^\alpha \delta_{j-1}^\beta $ for $i<j$ and for all $(\alpha,\beta)\in \{0,1\}^2$. If $p>q\geq 0$, then the set of morphisms $\square([p],[q])$ is empty. If $p = q$, then the set $\square([p],[p])$ is the singleton $\{\id_{[p]}\}$. For $0\leq p \leq q$, all maps of $\square$ from $[p]$ to $[q]$ are one-to-one. A good reference for presheaves is \cite{MR1300636}.

\bd \cite{Brown_cube} The category of presheaves over $\square$, denoted by $\square^{op}\set$, is called the category of {\rm precubical sets}.  A precubical set $K$ consists of a family of sets $(K_n)_{n \geq 0}$ and of set maps $\de_i^\alpha:K_n \rightarrow K_{n-1}$ with $1\leq i \leq n$ and $\alpha\in\{0,1\}$ satisfying the cubical relations $\de_i^\alpha\de_j^\beta = \de_{j-1}^\beta \de_i^\alpha$ for any $\alpha,\beta\in \{0,1\}$ and for $i<j$. An element of $K_n$ is called a {\rm $n$-cube}. Let $\dim(c)=n$ if $c\in K_n$. An element of $K_0$ is also called a {\rm vertex} of $K$. \ed

Let $K$ be a precubical set. There exists a functor $\square(K):(\square\ddownarrow K) \to \square^{op}\set$ which takes the map of precubical sets $\square[n]\to K$ to $\square[n]$. It is a general property of presheaves that $K = \liminj \square(K)$, and the latter colimit is denoted by $\liminj_{\square[n]\rightarrow K} \square[n]$. Let \[K_{\leq n} = \liminj_{\substack{\square[p]\to K\\p\leq n}} \square[p].\]  Let $\square[n]:=\square(-,[n])$. The boundary of $\square[n]$ is the precubical set $\square[n]_{\leq n-1}$ also denoted by $\de \square[n]$. In particular, one has $\de \square[0] = \varnothing$. 

\bd
A {\rm cocubical object} of a category $\C$ is a functor $\square\to \C$.
\ed

\begin{nota}
	Let $\C$ be a cocomplete category. Let $X:\square\to \C$ be a cocubical object of $\C$. Let \[\widehat{X}(K)=\liminj_{\square[n]\to K} X([n]).\]
\end{nota}

\bp \cite[Proposition~2.3.2]{realization} \label{eq-cat}
Let $\C$ be a cocomplete category. The mapping $X\mapsto \widehat{X}$ induces an equivalence of categories between the category of cocubical objects of $\C$ and the colimit-preserving preserving functors from $\square^{op}\set$ to $\C$.
\ep

Definition~\ref{def-rea-flow} is new. Moreover, only q-realization functors are implicitly studied in \cite{ccsprecub} because the h-model structure and the m-model structure of flows were not yet known: they are introduced 13 years later in \cite{QHMmodel}.

\bd \label{def-rea-flow} Let $r\in \{q,m,h\}$. A functor $F:\square^{op}\set \to \dtop$ is a {\rm r-realization functor} if it satisfies the following properties:
\begin{itemize}[leftmargin=*]
	\item $F$ is colimit-preserving.
	\item For all $n\geq 0$, the map of flows $F(\de\square[n])\to F(\square[n])$ is an r-cofibration.
	\item There is an objectwise weak equivalence of cocubical flows $F(\square[*])\to \{0<1\}^*$ in the r-model structure of $\dtop$.
\end{itemize}
\ed

\bp \label{elm}
Let $r\in \{q,m,h\}$. Let $F:\square^{op}\set \to \dtop$ be an r-realization functor. Then for all precubical sets $K$, the flow $F(K)$ is r-cofibrant and there is a natural bijection $K_0\iso F(K)^0$.
\ep

\bpf 
Let $K$ be a precubical set. Then the canonical map $\varnothing\to K$ is a transfinite composition of pushouts of the maps $\de\square[n]\to \square[n]$ for $n\geq 0$. Consequently, the canonical map $\varnothing\to F(K)$ is a transfinite composition of pushouts of the maps $F(\de\square[n])\to F(\square[n])$ for $n\geq 0$. It implies that $F(K)$ is r-cofibrant. From the objectwise weak equivalence of cocubical flows $F(\square[*])\to \{0<1\}^*$, we deduce the objectwise bijection of cocubical sets $F(\square[*])^0\iso \{0,1\}^* \iso \square[*]_0$. We obtain the natural bijection $F(K)^0 \iso K_0$.
\epf

\bth \label{homotopy-natural}
Let $r\in \{q,m,h\}$. Consider two r-realization functors \[F_1,F_2:\square^{op}\set \longrightarrow \dtop.\] Then there exists a natural transformation $\mu:F_1 \Rightarrow F_2$ such that there is a commutative diagram of cocubical flows 
\[
\xymatrix@C=2em@R=2em
{
F_1(\square[*]) \fd{}\fr{\mu_{\square[*]}} & F_2(\square[*]) \fd{} \\
{\{0<1\}^*} \ar@{=}[r] & {\{0<1\}^*}
}
\]
and such that for all precubical sets $K$, the map $\mu_K:F_1(K)\to F_2(K)$ natural with respect to $K$ is a weak equivalence of the r-model structure of $\dtop$. Moreover, for all $(\alpha,\beta)\in K_0\p K_0$, the natural map $\P_{\alpha,\beta}F_1(K)\stackrel{\simeq}\longrightarrow \P_{\alpha,\beta}F_2(K)$ is a homotopy equivalence of spaces.
\eth

\bpf
The maps of cocubical flows $F_i(\square[*])\to \{0<1\}^*$ for $i=1,2$ are objectwise fibrations since $\P_{\alpha,\beta}\{0<1\}^*$ is empty or equal to a singleton and because all topological spaces are fibrant. Consequently, they are objectwise trivial fibrations by definition of an r-realization functor. By \cite[Theorem~2.3.3]{realization}, there exists a natural transformation $\mu:F_1 \Rightarrow F_2$ such that there is the commutative diagram of cocubical ﬂows depicted in the statement of the theorem and such that, for all precubical sets $K$, the natural map $\mu_K:F_1(K) \Rightarrow F_2(K)$ is a simplicial homotopy equivalence, and therefore a weak equivalence of the r-model structure by \cite[Proposition~9.5.16]{ref_model2}, between two r-cofibrant flows. If $r=h$, then the natural map $\P_{\alpha,\beta}F_1(K)\to \P_{\alpha,\beta}F_2(K)$ is a homotopy equivalence of spaces by definition of the weak equivalences of the h-model structure of flows. If $r=q$, then the flows $F_1(K)$ and $F_2(K)$ are q-cofibrant by Proposition~\ref{elm}. Therefore, the spaces $\P_{\alpha,\beta}F_1(K)$ and $\P_{\alpha,\beta}F_2(K)$ are q-cofibrant by \cite[Theorem~5.7]{leftproperflow}. Using the Whitehead theorem \cite[Theorem~7.5.10]{ref_model2}, we deduce that the natural map $\P_{\alpha,\beta}F_1(K)\to \P_{\alpha,\beta}F_2(K)$ is a homotopy equivalence of spaces. It remains the case $r=m$. The flows $F_1(K)$ and $F_2(K)$ are m-cofibrant by Proposition~\ref{elm}. We deduce that the spaces $\P_{\alpha,\beta}F_1(K)$ and $\P_{\alpha,\beta}F_2(K)$ are m-cofibrant by \cite[Theorem~8.7]{QHMmodel}. By \cite[Corollary~3.4]{mixed-cole}, we deduce that the weak homotopy equivalence $\P_{\alpha,\beta}F_1(K)\to \P_{\alpha,\beta}F_2(K)$ is a homotopy equivalence of spaces as well.
\epf

\bth \label{exists}
There exists a q-realization functor $|\!-\!|_q:\square^{op}\set \to \dtop$.
\eth

\bpf
Let $(-)^{cof}$ be a q-cofibrant replacement functor of $\dtop$. Let 
\[{|K|_q := \liminj_{\square[n]\rightarrow
		K} (\{0<1\}^n)^{cof}}. \]
It is a q-realization functor by \cite[Proposition~7.4]{ccsprecub}. 
\epf

\begin{rem}
	The functor \[|\!-\!|_{bad}:K\mapsto \liminj_{\square[n]\rightarrow
		K} \{0<1\}^n\] is not a q-realization functor since the map of flows $|\de\square[2]|_{bad}\to |\square[2]|_{bad}$ is not a q-cofibration of flows.   \end{rem}

Moreover, there is the isomorphism $|\de\square[n]|_{bad}\iso |\square[n]|_{bad}$ for all $n\geq 3$ by \cite[Theorem~7.1]{ccsprecub}, which is not the expected behavior for a realization functor. 

\bp \label{qTOmTOh}
Every q-realization functor is an m-realization functor. Every m-realization functor is an h-realization functor.
\ep

\bpf
Every q-realization functor is an m-realization functor because every q-cofibration of flows is an m-cofibration of flows by \cite[Proposition~7.8]{QHMmodel} and because the weak equivalences are the same in the two model structures. Let $F:\square^{op}\set \to \dtop$ be an m-realization functor. Then for all $n\geq 0$, the map of flows $F(\de\square[n])\to F(\square[n])$ is an m-cofibration, and therefore an h-cofibration by \cite[Proposition~3.6]{mixed-cole}. The map of flows $F(\square[n])\to \{0<1\}^n$ is a weak equivalence of the m-model structure for all $n\geq 0$. Since $F(\square[n])$ is m-cofibrant by Proposition~\ref{elm}, there exists by \cite[Corollary~3.7]{mixed-cole} a q-cofibrant flow $C_n$ and a weak equivalence of the h-model structure of flows $C_n\to F(\square[n])$ for all $n\geq 0$. By \cite[Theorem~5.7]{leftproperflow}, for all $\alpha,\beta\in C_n^0 = \{0,1\}^n$, the topological space $\P_{\alpha,\beta}C_n$ is q-cofibrant. It means that for all $\alpha,\beta\in \{0,1\}^n$, the space $\P_{\alpha,\beta}F(\square[n])$ is homotopy equivalent to a q-cofibrant space, which means that $\P_{\alpha,\beta}F(\square[n])$ is m-cofibrant. Thus the map $\P_{\alpha,\beta} F(\square[n])\to \P_{\alpha,\beta}\{0<1\}^n$ is for all $\alpha,\beta\in \{0,1\}^n$ a weak homotopy equivalence between m-cofibrant spaces, and therefore a homotopy equivalence by \cite[Corollary~3.4]{mixed-cole}. In other terms, the map of flows $F(\square[n])\to \{0<1\}^n$ is a weak equivalence of the h-model structure of flows for all $n\geq 0$. We have proved that $F$ is an h-realization functor.
\epf

The drawback of the construction of Theorem~\ref{exists} is that it depends on the non-canonical choice of a q-cofibrant replacement. It is one of the purpose of the paper to fix this issue. The following theorem is not used in the sequel. It helps the reader to understand the geometric contents of a q-realization functor.

\bth \label{rea-func} \cite[Theorem~4.2.4 and Theorem~4.2.6]{realization}
For all $n\geq 1$, there is a homotopy pushout diagram of flows for the q-model structure
\[
\xymatrix@C=3em@R=3em{
	\glob(\mathbf{S}^{n-2}) \fr{\substack{0\mapsto 0_n\\1\mapsto 1_n}} \fd{}& |\de\square[n]|_q\fd{}\\
	\glob(\mathbf{D}^{n-1}) \fr{} & \hococartesien{|\square[n]|_q}.}
\]
There exists a q-realization functor such that the pushout diagram above is moreover strict.
\eth

\section{Natural \mins{d}-paths}
\label{natural_sec}

We want to use the notion of natural $d$-path introduced by Raussen in \cite[Definition~2.14]{MR2521708} to build the natural realization functor from precubical sets to flows. 

This new realization functor is natural in the sense that it uses natural $d$-paths, and also natural in the sense that it is more canonical than the q-realization functor of Theorem~\ref{exists}. Indeed, the latter depends on the non-canonical choice of a q-cofibrant replacement functor for the category of flows. The new one is independent of such a non-canonical choice.

\begin{nota}
	Let $\delta_i^\alpha : [0,1]^{n-1} \rightarrow [0,1]^n$ be the continuous map defined for $1\leq i\leq n$ and $\alpha \in \{0,1\}$ by $\delta_i^\alpha(\epsilon_1, \dots, \epsilon_{n-1}) = (\epsilon_1,\dots, \epsilon_{i-1}, \alpha, \epsilon_i, \dots, \epsilon_{n-1})$. By convention, let $[0,1]^0=\{()\}$. We obtain a cocubical topological space $[0,1]^*$ and the associated colimit-preserving functor from precubical sets to topological spaces is denoted by 
	\[
	|K|_{geom} = \liminj_{\square[n]\to K} [0,1]^n.
	\]
\end{nota}

The topological space $|K|_{geom}$ is a CW-complex, and therefore it is Hausdorff. Every point of $|K|_{geom}$ admits a unique presentation $[c;x]=|c|_{geom}(x)$ where $c$ is a cube of $K$ and such that $x\in ]0,1[^{\dim(c)}$. A point of $|K|_{geom}$ may belong to several cubes and therefore admits several presentations $[c;x]$ with $x\in [0,1]^{\dim (c)}$.

\bd
Let $U$ be a topological space. A {\rm (Moore) path} in $U$ consists of a continuous map $[0,\ell]\to U$ with $\ell>0$. The real number $\ell>0$ is called the {\rm length} of the path.
\ed

\bd
Let $\gamma_1:[0,\ell_1]\to U$ and $\gamma_2:[0,\ell_2]\to U$ be two paths in a topological space $U$ such that $\gamma_1(\ell_1)=\gamma_2(0)$. The {\rm Moore composition} $\gamma_1*\gamma_2:[0,\ell_1+\ell_2]\to U$ is the Moore path defined by 
\[
(\gamma_1*\gamma_2)(t)=
\begin{cases}
\gamma_1(t) & \hbox{ for } t\in [0,\ell_1]\\
\gamma_2(t-\ell_1) &\hbox{ for }t\in [\ell_1,\ell_1+\ell_2].
\end{cases}
\]
The Moore composition of Moore paths is strictly associative.
\ed

\bd  Let $n\geq 1$. A {\rm (nonconstant) tame $d$-path} of $|\square[n]|_{geom} = [0,1]^n$ is a nonconstant continuous map $\gamma:[0,\ell]\to [0,1]^n$ with $\ell>0$ such that $\gamma(0),\gamma(\ell)\in \{0,1\}^n$ and such that $\gamma$ is nondecreasing with respect to each axis of coordinates. 
\ed

\bd \label{dpath}
Let $K$ be a precubical set. A {\rm nonconstant tame $d$-path} of $K$ is a path $[0,\ell] \to |K|_{geom}$ which is the Moore composition $\gamma_1 * \dots *\gamma_n$ of nonconstant tame $d$-paths in cubes of $|K|_{geom}$. $\gamma(0)\in K_0$ is called the {\rm initial state} of $\gamma$ and $\gamma(\ell)\in K_0$ is called the {\rm final state} of $\gamma$. 
\ed

\begin{rem} \label{warning}
	All $d$-paths are tame and nonconstant in this paper. In particular, they start and end at a vertex of $K$. The adjective tame is added everywhere. The adjective nonconstant is often omitted.
\end{rem}

\begin{nota} \label{presentation}
	With the notations of Definition~\ref{dpath}. A tame $d$-path $\gamma:[0,\ell] \to |K|_{geom}$ can be written 
	$
	\gamma = {}^0[c_1;\gamma_1]\stackrel{t_1}*\dots \stackrel{t_{n-1}}*[c_n;\gamma_n]^{t_n}
	$
	or
	$
	\gamma = [c_1;\gamma_1]*\dots *[c_n;\gamma_n]
	$
	with $0=t_0<  t_1 < \dots < t_n=\ell$ such that for all $1\leq i\leq n$ and $t\in [t_{i-1},t_i]$, $\gamma(t)=[c_i;\gamma_i(t)]$ with $\dim(c_i)\geq 1$ and such that $\gamma(t_i)\in K_0$ for $0\leq i \leq n$. The sequence $(c_1,\dots,c_n)$ is called a {\rm carrier} of $\gamma$. The notation $\carrier(\gamma)$ means that a carrier of $\gamma$ is chosen: it is not unique. 
\end{nota}

The adjective \textit{tame} corresponds to the condition $\gamma(t_i)\in K_0$ for $0\leq i \leq n$. An important feature shared by all $d$-paths (tame or not) of a precubical set $K$ is that they have a well-defined $L_1$-arc length \cite[Section~2.2.1]{MR2521708} \cite[Section~2.2]{Raussen2012}. Intuitively, the natural $d$-paths are the $d$-paths whose speed corresponds to the $L_1$-arc length. We give an explicit definition of a tame natural $d$-path which is sufficient for this paper by starting from the tame $d$-paths in the topological $n$-cube $[0,1]^n$. It is equivalent to Raussen's definition of nonconstant tame natural $d$-path.

\bd 
Let $n\geq 1$. A {\rm tame natural $d$-path} of the topological $n$-cube $[0,1]^n$ is a $d$-path $\gamma=(\gamma_1,\dots,\gamma_n):[0,n]\to [0,1]^n$ such that for all $t\in [0,n]$, one has $t=\gamma_1(t)+\dots +\gamma_n(t)$. The set of tame natural $d$-paths in $[0,1]^n$ is denoted by $N_n$. It is equipped with the compact-open topology. 
\ed

\bd A tame $d$-path $\gamma$ of a precubical set $K$ is {\rm natural} if it can be written $\gamma = [c_1;\gamma_1]*\dots *[c_n;\gamma_n]$ such that each $\gamma_i$ is a tame natural $d$-path in the cube $c_i$ for all $i\in \{1,\dots,n\}$. 
\ed

\bp \label{Delta}
Let $n\geq 1$. The topological space $N_n$ is $\Delta$-generated and $\Delta$-Hausdorff. It is metrizable, contractible, compact and sequentially compact.
\ep

\bpf 
The compact-open topology is metrizable with the distance of the uniform convergence by \cite[Proposition~A.13]{MR1867354}. Therefore it is first countable. Consider a ball $B(\gamma,\epsilon)$ for this metric. Let $\gamma'\in B(\gamma,\epsilon)$. Then each convex combination $(1-u)\gamma+u\gamma'$ is a tame natural $d$-path since $(1-u)t+ut=t$ and for all $t\in [0,n]$ and all $i\in \{1,\dots,n\}$, one has 
\[
|((1-u)\gamma_i+u\gamma'_i)(t) - \gamma_i(t)| = u|\gamma'_i(t) - \gamma_i(t)| < u\epsilon \leq \epsilon.
\]
It means that the space $N_n$ is locally path-connected. By \cite[Proposition~3.11]{MR3270173}, it is $\Delta$-generated, and also $\Delta$-Hausdorff, being metrizable. It is contractible since there is a homotopy $H:[0,1] \p N_n\to N_n$ between the identity of $N_n$ and the constant map taking each tame natural $d$-path to the tame natural $d$-path $\delta:t\mapsto (t/n,t/n,\dots,t/n)$ given by the convex combination $H(u,\gamma) = u\delta +(1-u)\gamma$. It is compact by \cite[Proposition~9.5]{MR4070250} applied to the sequence $\underline{n}=(n)$. We want to give a different argument which does not use Lipschitz maps on metric spaces. Let $(\gamma^k)_{k\geq 0}=(\gamma_1^k,\dots,\gamma_n^k)_{k\geq 0}$ be a sequence of $N_n$. By a Cantor diagonalization argument, one can suppose that the sequence $(\gamma^k(r))_{k\geq 0}$ of $[0,1]^n$ converges to $(\gamma_1^\infty(r),\dots,\gamma_n^\infty(r))$ for all $r\in \mathbb{Q}\cap [0,n]$. Let $\gamma_i^-(x)=\sup\{\gamma_i^\infty(r)\mid r\in \mathbb{Q}\cap [0,x]\}$ and $\gamma_i^+(x)=\inf\{\gamma_i^\infty(r)\mid r\in \mathbb{Q}\cap [x,n]\}$. Then, by density of $\mathbb{Q}$, for all $x\in [0,n]$, one has $(\gamma_1^+(x)-\gamma_1^-(x)) + \dots + (\gamma_n^+(x)-\gamma_n^-(x)) = 0$. Thus, for all $x\in [0,n]$ and for all $1\leq i\leq n$, since $\gamma_i^+(x)-\gamma_i^-(x)\geq 0$, we deduce that $\gamma_i^+(x)=\gamma_i^-(x)$. It means that $\gamma_i^+=\gamma_i^-:[0,n]\to [0,1]$ is continuous for all $i\in \{1,\dots,n\}$. Consequently, each sequence $(\gamma_i^k)_{k\geq 0}$ converges pointwise for $1\leq i\leq n$. By the second Dini theorem, the convergence is uniform. Using \cite[Lemma~6.10]{Moore2}, we deduce that $(\gamma^k)_{k\geq 0}$ has a convergent subsequence. We deduce that $N_n$ is sequentially compact, hence compact, being metrizable.
\epf

\begin{nota}
	Let $\underline{x}=(x_1,\dots,x_n)$ and $\underline{x}'=(x'_1,\dots,x'_n)$ be two elements of $[0,1]^n$. Let
	\[
	d_\infty(\underline{x},\underline{x}') = \max_{1\leq i \leq n} |x_i-x'_i|.
	\]
\end{nota}

\bd
Let $n\geq 2$. Let $\mathcal{V}_n=\{0,1\}^n \backslash \{0_n,1_n\}$. Consider the continuous map $\phi:N_n\to [0,1]$ defined by $\phi(\gamma) = \min_{(t,v)\in [0,n]\p \mathcal{V}_n} d_\infty(\gamma(t),v)$. Let $\de N_n = \phi^{-1}(0)$ equipped with the relative topology.
\ed

\begin{nota}
	Let $\de N_0=N_0=\de N_1=\varnothing$. 
\end{nota}

There is the proposition: 

\bp \label{dDelta}
Let $n\geq 2$. The underlying set of $\de N_n$ is exactly the set of Moore compositions of tame natural $d$-paths in subcubes of $[0,1]^n$. For every $\gamma\in \de N_n$, $\gamma([0,n])$ is included in the boundary of $[0,1]^n$. The topology of $\de N_n$ is $\Delta$-generated and $\Delta$-Hausdorff. It is metrizable, compact and sequentially compact.
\ep

\bpf The set $\de N_n$ is exactly the set of tame natural $d$-paths in $[0,1]^n$ whose image intersects $\mathcal{V}_n$. Let $\gamma=(\gamma_1,\dots,\gamma_n)\in \de N_n$ and let $t_0\in ]0,n[$ such that $\gamma(t_0)=(\epsilon_1,\dots,\epsilon_n)\in \mathcal{V}_n$. Since $\gamma$ is natural, one has $t_0=\epsilon_1+\dots+\epsilon_n$ which is therefore an integer between $1$ and $n-1$. Then $\gamma=\gamma^a*\gamma^b$ with $\gamma^a(0)=0_n$, $\gamma(t_0)=\gamma^a(t_0)=\gamma^b(0)\in \mathcal{V}_n$ and $\gamma^b(n-t_0)=1_n$. Therefore, for all $t\in [0,t_0]$, $t=\gamma_1(t) + \dots + \gamma_n(t)= \gamma^a_1(t) + \dots + \gamma^a_n(t)$. Let $J=\{j\in \{1,\dots,n\}\mid \epsilon_j=0\}$. Since the paths are nondecreasing with respect to each axis of coordinates, it implies that $\gamma^a_j(t)=0$ for all $j\in J$. Thus for all $t\in [0,t_0]$, $t=\sum_{j\notin J}\gamma^a_j(t)$. It means that $\gamma^a$ is a natural path in the subcube from $0_n$ to $\gamma^a(t_0)$. For all $t\in [t_0,n]$, on has $t=\gamma_1(t) + \dots + \gamma_n(t)=\gamma^b_1(t-t_0) + \dots + \gamma^b_n(t-t_0)$, the first equality since $\gamma$ is natural, the second equality by definition of the Moore composition of paths. We deduce that $t-t_0=(\gamma^b_1(t-t_0) -\epsilon_1) + \dots + (\gamma^b_n(t-t_0)-\epsilon_n)$ for all $t\in [t_0,n]$. If for some $i\in \{1,\dots,n\}$, $\epsilon_i=1$, then $\gamma^b_i=1$ since the paths are nondecreasing with respect to each axis of coordinates.  We obtain $t-t_0=\sum_{j\in J} \gamma^b_j(t-t_0)$ for all $t\in [t_0,n]$. It means that $\gamma^b$ is a natural path in the subcube going from $\gamma(t_0)$ to $1_n$. We deduce that the underlying set of $\de N_n$ is exactly the set of Moore compositions of tame natural $d$-paths in subcubes of $[0,1]^n$. The second assertion is a consequence of this fact. Consider $\gamma \in \de N_n$. There exists $t_0\in ]0,n[$ such that $\gamma(t_0)\in \mathcal{V}_n$. Since $\mathcal{V}_n$ is discrete, there exists an open $U$ of $[0,1]^n$ such that $U\cap \mathcal{V}_n = \{\gamma(t_0)\}$. Then $W(\{t_0\},U) = \{\gamma'\in \de N_n\mid \gamma'(t_0)\in U\}$ is an open subset of $\de N_n$ for the compact-open topology. The latter being metrizable, take a ball $B(\gamma,\epsilon) \subset W(\{t_0\},U)$ and repeat the reasoning of Proposition~\ref{Delta}: we deduce that $\de N_n$ is locally path-connected as well, hence $\Delta$-generated (and also $\Delta$-Hausdorff, being metrizable) by \cite[Proposition~3.11]{MR3270173}.
\epf

To summarize, $\de N_n$ is a closed subset of $N_n$ which remains $\Delta$-generated and also $\Delta$-Hausdorff when equipped with the relative topology. Both $\de N_n$ and $N_n$ are equipped with the compact-open topology and are metrizable, compact, and sequentially compact. 

The presentation chosen for $\de N_n$ and $N_n$ is due to the fact that \cite[Proposition~10.2]{MR4070250} is used in the proof of Proposition~\ref{deISde} and that \cite[Proposition~10.3]{MR4070250} is used in the proof of Proposition~\ref{h-cof}. However, \cite{MR4070250} uses the compact-open topology. It turns out that the $\Delta$-kelleyfication functor does not preserve compactness and that it is a right Quillen equivalence, not a left Quillen equivalence. Since $\de N_n$ and $N_n$ equipped with the compact-open topology are metrizable for all $n\geq 0$, they are k-spaces. An additional argument is necessary to prove that they are $\Delta$-generated as well to get rid of this issue.

\section{Natural realization from precubical sets to flows}
\label{nat_rea_sec}

We define a flow $|\square[n]|_{nat}$ for $n\geq 0$ called \textit{the natural $n$-cube} as follows. The set of states is $\{0,1\}^n$. Let $n\geq 1$ and $\alpha,\beta\in\{0,1\}^n$. Let $\alpha=(\alpha_1,\dots,\alpha_n)$ and $\beta=(\beta_1,\dots,\beta_n)$. Assume that $\alpha<\beta$ in the product order $\{0<1\}^n$. Let $I=\{i\in \{1,\dots,n\}\mid \alpha_i \neq \beta_i\}$. By hypothesis, $I$ is nonempty. Let $m$ be the cardinality of $I$. Then $\alpha$ ($\beta$ resp.) is the initial (final resp.) state of a $m$-subcube of $\square[n]$. Then let $\P_{\alpha,\beta}|\square[n]|_{nat} = N_m$ viewed as the space of tame natural $d$-paths in the $m$-subcube from $\alpha$ to $\beta$. Assume that $\alpha\geq \beta$. Let $\P_{\alpha,\beta}|\square[n]|_{nat} = \varnothing$. The composition law is defined by the Moore composition of tame natural $d$-paths, which is still a tame natural $d$-path.

\bp \label{incl}
Let $\phi:[n]\to [n+1]$ be a map of the small category $\square$. Then the continuous map $\P_{0_{n},1_{n}}|\square[n]|_{nat}\to \P_{\phi(0_{n}),\phi(1_{n})}|\square[n+1]|_{nat}$ induced by $\phi$ is the identity of $N_n$.
\ep

\bpf
It is a straightforward consequence of the definitions.
\epf

\begin{cor}
	We obtain a well-defined cocubical flow $|\square[*]|_{nat}$. 
\end{cor}

\bpf
Consider an algebraic relation $\phi_1\phi_2=\psi_1\psi_2:[n] \to [n+2]$ in the small category $\square$. Consider the diagram of topological spaces 
\[
\xymatrix@C=2em@R=2em
{
\P_{0_{n},1_{n}}|\square[n]|_{nat}\ar@{=}[d] \fr{} & \P_{\phi_2(0_{n}),\phi_2(1_{n})}|\square[n+1]|_{nat} \fr{} & \P_{\phi_1\phi_2(0_{n}),\phi_1\phi_2(1_{n})}|\square[n+2]|_{nat} \ar@{=}[d]\\
\P_{0_{n},1_{n}}|\square[n]|_{nat}\fr{} & \P_{\psi_2(0_{n}),\psi_2(1_{n})}|\square[n+1]|_{nat} \fr{} & \P_{\psi_1\psi_2(0_{n}),\psi_1\psi_2(1_{n})}|\square[n+2]|_{nat}
}
\] 
By Proposition~\ref{incl} and by definition of $|\!-\!|_{nat}$, the two horizontal composite maps are equal to the identity of $N_n$. It means that the diagram is commutative and that the cocubical relations are satisfied.
\epf

Using Proposition~\ref{eq-cat}, we obtain:

\bd \label{natural-rea}
Let $K$ be a precubical set. Consider the colimit-preserving functor 
\[
|K|_{nat} = \liminj_{\square[n]\rightarrow K} |\square[n]|_{nat}.
\]
It is called the {\rm natural realization of $K$} as a flow.
\ed

\bp \label{deISde}
Let $n\geq 0$. There is a homeomorphism $\de N_n \iso \P_{0_{n},1_{n}}|\de\square[n]|_{nat}$.
\ep

\bpf
Using \cite[Proposition~10.2]{MR4070250} applied with the sequence $(n)$, we deduce that this map is a homeomorphism: the idea of the proof is that there is a continuous bijection from $\P_{0_{n},1_{n}}|\de\square[n]|_{nat}$ to $\de N_n$ and that both $\P_{0_{n},1_{n}}|\de\square[n]|_{nat}$ and $\de N_n$ are compact.
\epf

\bp \label{h-cof}
The continuous map $\de N_n \subset N_n$ is an h-cofibration of spaces for all $n\geq 0$. 
\ep

\bpf
Using \cite[Proposition~10.3]{MR4070250} applied with the sequence $(n)$, we deduce that this map is a strong neighborhood deformation retract, i.e. an h-cofibration by \cite[Theorem~2]{vstrom1}.
\epf

\begin{cor} \label{h-cof-glob}
	The map of flows $\glob(\de N_n) \subset \glob(N_n)$ is an h-cofibration of flows for all $n\geq 0$.
\end{cor}

\bpf
A map of flows of the form $\glob(U)\to \glob(V)$ satisfies the left lifting property with respect to a map of flows $f:X\to Y$ if and only if the map $U\to V$ satisfies the left lifting property with respect to all maps $\P_{\alpha,\beta}X\to \P_{f(\alpha),f(\beta)}Y$ for all $(\alpha,\beta)\in X^0\p X^0$. Using the characterization of the trivial h-fibrations of flows (see Theorem~\ref{three}) and Proposition~\ref{h-cof}, we deduce that the map $\glob(\de N_n)\to \glob(N_n)$ is an h-cofibration of flows.
\epf

\bp \label{cof-iteration}
For all $n\geq 0$, the map $|\de \square[n]|_{nat} \to |\square[n]|_{nat}$ is an h-cofibration of flows. 
\ep

\bpf
From the homeomorphism of Proposition~\ref{deISde} and by definition of $N_n$, we deduce that the commutative diagram of spaces
\[
\xymatrix@C=3em@R=3em
{
	\de N_n \fr{\iso} \fd{} &  \P_{0_{n},1_{n}}|\de \square[n]|_{nat} \fd{} \\
	N_n\fr{\iso} & \cocartesien {\P_{0_{n},1_{n}}|\square[n]|_{nat}}
}
\]
is a pushout diagram of spaces. The top homeomorphism yields a map of flows \[\glob(\de N_n)\longrightarrow |\de \square[n]|_{nat}\] taking $0$ to $0_n$ and $1$ to $1_n$ for all $n\geq 0$. We obtain the pushout diagram of flows
\[
\xymatrix@C=3em@R=3em
{
	\glob(\de N_n) \fr{} \fd{} & |\de \square[n]|_{nat} \fd{} \\
	\glob(N_n) \fr{} & \cocartesien {|\square[n]|_{nat}}
}
\]
Using Corollary~\ref{h-cof-glob}, we deduce that the map $|\de \square[n]|_{nat} \to |\square[n]|_{nat}$ is an h-cofibration of flows for all $n\geq 0$.
\epf

\bth \label{natural-equal}
There exists a natural transformation $\mu:|\!-\!|_q\Rightarrow |\!-\!|_{nat}$ such that for all precubical sets $K$, the natural map $\mu_K:|K|_q \to |K|_{nat}$ induces a bijection on states and a homotopy equivalence $\P_{\alpha,\beta}|K|_q \simeq \P_{\alpha,\beta}|K|_{nat}$ for all $\alpha,\beta\in K_0$.
\eth

\bpf
The map $|\square[*]|_{nat}\to \{0<1\}^*$ is an objectwise weak equivalence for the h-model structure of $\dtop$ since all spaces $N_n$ for $n\geq 1$ are contractible by Proposition~\ref{Delta}. By Proposition~\ref{cof-iteration}, the natural realization functor is then an h-realization functor. Since $|\!-\!|_q$ is also an h-realization functor by Proposition~\ref{qTOmTOh}, the proof is complete thanks to Theorem~\ref{homotopy-natural}.
\epf

The statement of Theorem~\ref{homotopy-natural} being symmetric, there is also a natural transformation $\nu:|\!-\!|_{nat} \Rightarrow |\!-\!|_q$ such that, for all precubical sets $K$, the natural map $\nu_K:|K|_{nat} \to |K|_{q}$ induces a bijection on states and a homotopy equivalence $\P_{\alpha,\beta}|K|_{nat} \simeq \P_{\alpha,\beta}|K|_{q}$ for all $\alpha,\beta\in K_0$. This statement is less intuitive because, morally speaking, the natural realization contains more execution paths than the q-realization.

Proposition~\ref{cof-iteration} means that the natural realization functor is an h-realization functor. In fact, it is possible to prove better. For all precubical sets $K$ and all $(\alpha,\beta)\in K_0\p K_0$, there is a homotopy equivalence $\P_{\alpha,\beta}|K|_{nat} \simeq \P_{\alpha,\beta}|K|_q$. Since $|K|_q$ is q-cofibrant, the space $\P_{\alpha,\beta}|K|_q$ is q-cofibrant by \cite[Theorem~5.7]{leftproperflow}. It means that the spaces of execution paths $\P_{\alpha,\beta}|K|_{nat}$ are m-cofibrant for all $(\alpha,\beta)\in K_0\p K_0$. This suggests that the natural realization $|K|_{nat}$ is an m-cofibrant flow. Indeed we have the following theorem:

\bth \label{m-cof}
The natural realization functor is an m-realization functor. For any precubical set $K$, the flow $|K|_{nat}$ is m-cofibrant. 
\eth

\bpf
The map $|\square[*]|_{nat}\to \{0<1\}^*$ is an objectwise weak equivalence for the h-model structure of $\dtop$, and therefore for the m-model structure of $\dtop$ as well. There is a homeomorphism $\de N_n\iso \P_{0_{n},1_{n}}|\de\square[n]|_{nat}$ (Proposition~\ref{deISde}) and a homotopy equivalence $\P_{0_{n},1_{n}}|\de\square[n]|_{nat} \simeq \P_{0_{n},1_{n}}|\de\square[n]|_q$ (Theorem~\ref{homotopy-natural}). Since $|\de\square[n]|_q$ is a q-cofibrant flow by Proposition~\ref{elm}, the space $\P_{0_{n},1_{n}}|\de\square[n]|_q$ is q-cofibrant by \cite[Theorem~5.7]{leftproperflow}. Moreover, $N_n$ is contractible by Proposition~\ref{Delta}, hence m-cofibrant. It implies that all maps $\de N_n\to N_n$ for $n\geq 0$ are h-cofibrations of spaces between m-cofibrant spaces \cite[Corollary~3.7]{mixed-cole}. By \cite[Corollary~3.12]{mixed-cole}, the maps $\de N_n\to N_n$ are therefore m-cofibrations of spaces for all $n\geq 0$. Thus, the map of flows $\glob(\de N_n)\to \glob(N_n)$ is an m-cofibration of flows for all $n\geq 0$ by the same argument as in the proof of Corollary~\ref{h-cof-glob}. Using the pushout diagram in the proof of Proposition~\ref{cof-iteration}, we deduce that the natural realization functor is an m-realization functor. By Proposition~\ref{elm}, we deduce that the flow $|K|_{nat}$ is m-cofibrant. 
\epf

\section{Natural realization and cube chains}
\label{cube_sec}

Cube chains are introduced in \cite[Definition~1.1]{MR3722069}. We use the presentation given in \cite[Section~7]{MR4070250} instead. Let $\seq(n)$ be the set of sequences of positive integers $\underline{n}=(n_1,\dots,n_p)$ with $n_1+\dots + n_p=n$. Let $\underline{n}=(n_1,\dots,n_p) \in \seq(n)$. Then $|\underline{n}|=n$ is the length of $\underline{n}$ and $\ell(\underline{n})=p$ is the number of elements of $\underline{n}$. Let $K$ be a precubical set and $A={a_1<\dots<a_k} \subset \{1,\dots,n\}$ and $\epsilon\in \{0,1\}$. The \textit{iterated face map} is defined by \[\de^\epsilon_A=\de^\epsilon_{a_1} \de^\epsilon_{a_2} \dots \de^\epsilon_{a_k}.\]

\bd 
Let $\underline{n}\in \seq(n)$. The {\rm $\underline{n}$-cube} is the precubical set 
\[
\square[\underline{n}] = \square[n_1] * \dots * \square[n_p]
\]
where the notation $*$ means that the final state $1_{n_i}$ of the precubical set $\square[n_i]$ is identified with the initial state $0_{n_{i+1}}$ of the precubical set $\square[n_{i+1}]$ for $1\leq i\leq p-1$.
\ed 

Let $K$ be a precubical set. Let $\alpha,\beta$ be two vertices of $K$. Let $n\geq 1$. The category $\Ch_{\alpha,\beta}(K,n)$ is defined as follows. The objects are the maps of precubical sets $\square[\underline{n}] \to K$ with $|\underline{n}|=n$ where the initial state of $\square[n_1]$ is mapped to $\alpha$ and the final state of $\square[n_p]$ is mapped to $\beta$. Let $A\sqcup B=\{1,\dots,m_1+m_2\}$ be a partition with the cardinal of $A$ equal to $m_1>0$ and the cardinal of $B$ equal to $m_2>0$. Let \[\phi_{A,B}:\square[m_1]*\square[m_2] \longrightarrow \square[m_1+m_2]\] be the unique map of precubical sets such that 
\begin{align*}
&\phi_{A,B}(\id_{[m_1]}) = \de^0_B(\id_{[m_1+m_2]}),\\
&\phi_{A,B}(\id_{[m_2]}) = \de^1_A(\id_{[m_1+m_2]}).
\end{align*}
For $i\in \{1,\dots,\ell(\underline{n})\}$ and a partition $A\sqcup B=\{1,\dots,n_i\}$, let \[\delta_{i,A,B}=\id_{\square[n_1]}*\dots*\id_{\square[n_{i-1}]}*\phi_{A,B}*\id_{\square[n_{i+1}]}*\dots*\id_{\square[n_{\ell(\underline{n})}]}.\] The morphisms are the commutative diagrams
\[
\xymatrix@C=3em@R=3em
{
	\square[\underline{n}_a] \fd{}\fr{a} & K \ar@{=}[d] \\
	\square[\underline{n}_b]  \fr{b}  & K
}
\] 
where the left vertical map is a composite of maps of precubical sets of the form $\delta_{i,A,B}$. 

From a precubical set $K$, we are going to define a flow $||K||$ as follows. The set of states is $K_0$. Consider the small diagram of spaces \[\D_{\alpha,\beta}(K,n):\Ch_{\alpha,\beta}(K,n)\longrightarrow \top\] defined by on objects by \[\D_{\alpha,\beta}(K,n)(\square[\underline{n}]\to K) = N_{n_1} \p \dots \p N_{n_p}\] and on morphisms by using the maps \[\P|\phi_{A,B}|_{nat}:\P(\square[m_1]*\square[m_2]) \to \P(\square[m_1+m_2])\] which induce maps $N_{m_1}\p N_{m_2}\to N_{m_1+m_2}$ given by the Moore composition of tame natural $d$-paths. The space of execution spaces $\P_{\alpha,\beta} ||K||$ is defined as follows:
\[
\P_{\alpha,\beta} ||K|| = \displaystyle\coprod_{n\geq 1} \liminj \D_{\alpha,\beta}(K,n).
\]
It is easy to see that the concatenation of tuples induces functors \[\D_{\alpha,\beta}(K,m_1)\p \D_{\beta,\gamma}(K,m_2) \to \D_{\alpha,\gamma}(K,m_1+m_2),\]
and, using \cite[Proposition~A.4]{leftproperflow}, continuous maps 
\[\liminj \D_{\alpha,\beta}(K,m_1)\p \liminj \D_{\beta,\gamma}(K,m_2) \to \liminj\D_{\alpha,\gamma}(K,m_1+m_2)\] for all $m_1,m_2\geq 1$. We obtain an associative composition map \[\P_{\alpha,\beta} ||K|| \p \P_{\beta,\gamma} ||K|| \to \P_{\alpha,\gamma} ||K||\] for all $(\alpha,\beta,\gamma)\in K_0\p K_0\p K_0$.

\bp \label{lem_prep}
There is an isomorphism of cocubical flows \[||\square[*]|| \iso |\square[*]|_{nat}.\]
\ep

\bpf 
At first, we prove the isomorphism of flows $||\square[n]|| \iso |\square[n]|_{nat}$ by induction on $n\geq 0$. The statement is obvious for $n=0$. Let $n\geq 1$ and $\alpha,\beta\in\{0,1\}^n$ with $\alpha<\beta$. Let $\alpha=(\alpha_1,\dots,\alpha_n)$ and $\beta=(\beta_1,\dots,\beta_n)$. Let \[I=\{i\in \{1,\dots,n\}\mid \alpha_i \neq \beta_i\}.\] By hypothesis, $I$ is nonempty. Let $m$ be the cardinality of $I$. Then $\alpha$ ($\beta$ resp.) is the initial (final resp.) state of a $m$-subcube $\underline{c}$ of $\square[n]$. We deduce that the category $\Ch_{\alpha,\beta}(\square[n],p)$ is empty for $p\neq m$ and that it has a terminal object $\underline{c}:\square[m]\to \square[n]$ for $p=m$ corresponding to the subcube from $\alpha$ to $\beta$. We deduce the homeomorphisms
\begin{align*}
\P_{\alpha,\beta} ||\square[n]|| &= \liminj_{\substack{\underline{n}=(n_1,\dots,n_p),\ell(\underline{n})=m\\\square[\underline{n}]\to \square[n] \in \Ch_{\alpha,\beta}(\square[n],m)}} N_{n_1} \p \dots \p N_{n_p} \\&\iso N_m \\&= \P_{\alpha,\beta} |\square[n]|_{nat},
\end{align*}
the first equality by definition of $||\square[n]||$, the homeomorphism because of the unique map $\underline{c}:\square[m]\to \square[n]$ which is the terminal object of $\Ch_{\alpha,\beta}(\square[n],m)$, and the last equality by Proposition~\ref{incl} applied to the map $\underline{c}:\square[m]\to \square[n]$. By Proposition~\ref{incl} again, the isomorphism $||\square[n]|| \iso |\square[n]|_{nat}$ is natural with respect to $[n]$. 
\epf

We do not know yet that the functor $||\!-\!||$ is colimit-preserving. An additional argument based on Proposition~\ref{incl} as well is necessary for proving Theorem~\ref{iso-cube}.

\bth \label{iso-cube}
There is a natural isomorphism of flows \[||K|| \iso |K|_{nat}\] for all precubical sets $K$.
\eth

\bpf
Let $\underline{n}=(n_1,\dots,n_p) \in \seq(n)$. Every map of precubical sets $\square[\underline{n}]\to K$ gives rise to a map of flows $|\square[\underline{n}]|_{nat}\to |K|_{nat}$, and therefore to a continuous map \[N_{n_1}\p \dots \p N_{n_p} \longrightarrow \P  |K|_{nat}.\] Let $\phi_{A,B}:\square[m_1]*\square[m_2] \to \square[m_1+m_2]$ as above. A composite map of precubical sets $\square[m_1]*\square[m_2] \to \square[m_1+m_2] \to K$ gives rise to the commutative diagram of flows 
\[
\xymatrix@C=3em@R=3em
{|\square[m_1]*\square[m_2]|_{nat} \fd{}\fr{} & |K|_{nat} \ar@{=}[d] \\
	|\square[m_1+m_2]|_{nat}\fr{} & |K|_{nat}}
\]
and therefore to the commutative diagram of spaces
\[
\xymatrix@C=3em@R=3em
{N_{m_1} \p N_{m_2} \fd{}\fr{} & \P |K|_{nat} \ar@{=}[d] \\
	N_{m_1 + m_2} \fr{} & \P |K|_{nat}}
\]
Consequently, we obtain a cocone \[(N_{n_1} \p \dots \p N_{n_p})_{\substack{\square[\underline{n}]\to K \\\in \Ch_{\alpha,\beta}(K,n)}} \stackrel{\bullet}\longrightarrow \P |K|_{nat}\] and then a map of flows $||K|| \to |K|_{nat}$ which is bijective on states. For each map of precubical sets $\square[n]\to K$, we obtain by Proposition~\ref{lem_prep} a map of flows \[|\square[n]|_{nat} \iso ||\square[n]|| \longrightarrow ||K||.\] Using Proposition~\ref{incl}, we obtain a cocone of flows 
\[
(|\square[n]|_{nat})_{\square[n]\to K} \stackrel{\bullet}\longrightarrow ||K||
\]
and therefore a map of flows $|K|_{nat} \to ||K||$ such that the composite map $|K|_{nat} \to ||K|| \to |K|_{nat}$ is the identity. Thus, the map $||K||\to |K|_{nat}$ is onto on execution paths and the map $|K|_{nat} \to ||K||$ is one-to-one on execution paths. Consider an execution path $\gamma$ of $||K||$. It belongs to a colimit and therefore has a representative $(\gamma_1,\dots,\gamma_p)$ in a space of the form $N_{n_1} \p \dots \p N_{n_p}$ corresponding to some map of precubical set $\square[\underline{n}]\to K$ with $\underline{n}=(n_1,\dots,n_p)$. By Proposition~\ref{incl}, there exists an execution path $\overline{\gamma_1}*\dots * \overline{\gamma_p}$ of $|K|_{nat}$ which is mapped to $(\gamma_1,\dots,\gamma_p)$ by the map of flows $|K|_{nat}\to ||K||$. The latter is therefore surjective on execution paths. And the proof is complete.
\epf

\section{Comparing execution paths and \mins{d}-paths}
\label{spatial_precubical_set}

The fact that the Moore composition of tame natural $d$-paths in the geometric realization of a precubical set is strictly associative entails the following definition.

\bd \label{rea_conc}
Let $K$ be a precubical set. The {\rm tame concrete realization} of $K$ is the flow $|K|_{tc}$ such that the set of states is $K_0$, the space of execution paths $\P_{\alpha,\beta} |K|_{tc}$ for $(\alpha,\beta)\in K_0\p K_0$ is the space of nonconstant tame natural $d$-paths from $\alpha$ to $\beta$ in $|K|_{geom}$ equipped with the $\Delta$-kelleyfication of the relative topology induced by the compact-open topology and the composition of execution paths is induced by the Moore composition.
\ed

This construction yields a well-defined functor $|\!-\!|_{tc}:\square^{op}\set\to \dtop$. By definition of $|\square[n]|_{nat}$ (see the beginning of Section~\ref{nat_rea_sec}), there is a natural isomorphism of flows $|\square[n]|_{nat}\iso |\square[n]|_{tc}$ for all $[n]\in \square$. Since the natural realization functor is colimit-preserving, the universal property of the colimit provides a natural map of flows $|K|_{nat}\to |K|_{tc}$. 

\bp \label{prespatial}
Let $K$ be a precubical set. The natural map of flows $|K|_{nat}\to |K|_{tc}$ is bijective on states. For all $(\alpha,\beta)\in K_0\p K_0$, the continuous map $\P_{\alpha,\beta}|K|_{nat}\to \P_{\alpha,\beta}|K|_{tc}$ is onto. The following assertions are equivalent. 
\begin{enumerate}
	\item For all $(\alpha,\beta)\in K_0\p K_0$, the continuous map $\P_{\alpha,\beta}|K|_{nat}\to \P_{\alpha,\beta}|K|_{tc}$ is one-to-one.
	\item For all $(\alpha,\beta)\in K_0\p K_0$, the continuous map $\P_{\alpha,\beta}|K|_{nat}\to \P_{\alpha,\beta}|K|_{tc}$ is bijective.
%	\item For all $(\alpha,\beta)\in K_0\p K_0$, the continuous map of \eqref{emb} is one-to-one.
\end{enumerate}
Finally, for all $(\alpha,\beta)\in K_0\p K_0$, the space of execution paths $\P_{\alpha,\beta}|K|_{tc}$ is Hausdorff.
\ep

\bpf
The first assertion is by definition of $|K|_{nat}$ and $|K|_{tc}$. Let $(\alpha,\beta)\in K_0\p K_0$. Consider a tame natural $d$-path $[c_1;\gamma_1]*\dots *[c_p;\gamma_p]$ of $\P_{\alpha,\beta}|K|_{tc}$ (cf. Notation~\ref{presentation}). Let $n=\sum_i \dim(c_i)$ and $\underline{n}=(\dim(c_1),\dots,\dim(c_p))$. The sequence of cubes $(c_1,\dots,c_p)$ gives rise to a map $\square[\underline{n}] \to K$. Then \[(\gamma_1,\dots,\gamma_p)\in \D_{\alpha,\beta}(K,n)(\square[\underline{n}]\to K) = N_{\dim(c_1)} \p \dots \p N_{\dim(c_p)}\] is the representative of an element of $\P_{\alpha,\beta}|K|_{nat} \iso \P_{\alpha,\beta}||K||$ which is taken to $[c_1;\gamma_1]*\dots *[c_n;\gamma_n]\in\P_{\alpha,\beta}|K|_{tc}$. It means that, for all $(\alpha,\beta)\in K_0\p K_0$, the continuous map $\P_{\alpha,\beta}|K|_{nat}\to \P_{\alpha,\beta}|K|_{tc}$ is onto. It implies $(1)\Leftrightarrow (2)$. From the sequence of one-to-one continuous maps 
\[
\P_{\alpha,\beta}|K|_{tc} \subset \coprod\limits_{n\geq 1}\ttop([0,n],|K|_{geom}) \subset \coprod\limits_{n\geq 1}\top_{co}([0,n],|K|_{geom}) \subset \coprod\limits_{n\geq 1}\prod_{[0,n]} |K|_{geom}
\]
for all $(\alpha,\beta)\in K_0\p K_0$, $\top_{co}([0,n],|K|_{geom})$ being the set $\top([0,n],|K|_{geom})$ equipped with the compact-open topology, and the product being equipped with the pointwise topology (i.e. the product topology), we deduce that the space of execution paths $\P_{\alpha,\beta}|K|_{tc}$ is Hausdorff for all $(\alpha,\beta)\in K_0\p K_0$, $|K|_{geom}$ being Hausdorff. 
\epf

\bd \label{spatial}
A precubical set $K$ is {\rm spatial} if for all $(\alpha,\beta)\in K_0\p K_0$, the continuous map $\P_{\alpha,\beta}|K|_{nat}\to \P_{\alpha,\beta}|K|_{tc}$ is bijective.
\ed

A characterization of spatial precubical sets as a small orthogonality class is postponed to Appendix~\ref{spatial_app}. We want to give some examples of spatial precubical sets before expounding the main theorems of the section.

\bd \cite[page~499]{MR3722069}
A precubical set $K$ is {\rm proper} if the map 
\[
\coprod_{n\geq 0} K_n \longrightarrow K_0 \p K_0
\]
which takes an $n$-cube $c$ of $K$ to $(\de^0_{\{1,\dots,n\}}c,\de^1_{\{1,\dots,n\}}c)$ is one-to-one.
\ed

For all $n\geq 0$, the precubical sets $\de\square[n]$ and $\square[n]$ are proper. The precubical sets associated to all PV-programs are proper. Every geometric precubical set in the sense of \cite[Definition~1.18]{zbMATH07226006} is proper. In particular, every non-positively curved precubical set in the sense of \cite[Definition~1.28]{zbMATH07226006} is proper since it is geometric by definition. Indeed, let $K$ be a geometric precubical set. Let $c_1,c_2$ be two cubes of $K$ such that $(\de^0_{\{1,\dots,\dim(c_1)\}}c_1,\de^1_{\{1,\dots,\dim(c_1)\}}c_1) = (\de^0_{\{1,\dots,\dim(c_2)\}}c_2,\de^1_{\{1,\dots,\dim(c_2)\}}c_2)$. Suppose that $\dim(c_1)=0$. Since $K$ has no self-intersection, $c_2$ is $0$-dimensional as well and $c_1=c_2$. Assume that $\dim(c_1)\geq 1$. Then $c_2$ cannot be $0$-dimensional because $K$ has no self-intersection. The cubes $c_1$ and $c_2$ have a maximal common face which is necessarily $c_1=c_2$. Thus $K$ is proper.

\bp \label{referee}
There are two strict inclusions 
\[
\{\hbox{proper precubical sets}\} \subset \{\hbox{spatial precubical sets}\} \subset \{\hbox{precubical sets}\}.
\]
\ep

\bpf
Let $K$ be a proper precubical set. Let $(\alpha,\beta)\in K_0\p K_0$. Let $\xi_1$ and $\xi_2$ be two execution paths of $\P_{\alpha,\beta}|K|_{nat}$. Suppose that $\xi_1$ and $\xi_2$ are taken to the same tame natural $d$-path \[\gamma = [c_1;\gamma_1]*\dots *[c_p;\gamma_p]\] in $|K|_{geom}$ with $p\geq 1$ and $0=t_0<  t_1 < \dots < t_p=\ell$ such that for all $1\leq i\leq p$ and $t\in [t_{i-1},t_i]$, $\gamma(t)=[c_i;\gamma_i(t)]$ with $\dim(c_i)\geq 1$ and such that $\gamma(t_i)\in K_0$ for $0\leq i \leq p$. Choose the presentation $[c_1;\gamma_1]*\dots *[c_p;\gamma_p]$ so that $\gamma([0,1])\cap K_0=\{\gamma(t_i)\mid 0\leq i \leq p\}$ and $\gamma(]t_{i-1},t_i[) \cap K_0 = \varnothing$ for $1\leq i \leq p$. Let $n=\sum_i \dim(c_i)$ and  $\underline{n}=(\dim(c_1),\dots,\dim(c_p))$. There exist maps of precubical sets $a_k: \square[\underline{n}] \to K$ for $k=1,2$ such that \[(\gamma_1,\dots,\gamma_p)\in \D_{\alpha,\beta}(K,n)(a_k)=N_{\dim(c_1)}\p \dots \p N_{\dim(c_p)}\] is identified to $\xi_k$ for $k=1,2$ in the colimit $\P_{\alpha,\beta}||K||$. Since $K$ is proper and since $\gamma(]t_{i-1},t_i[) \cap K_0 = \varnothing$ for $1\leq i \leq p$, we deduce that $a_1=a_2$. Therefore, $\xi_1=\xi_2$. It means that $K$ is spatial. Consider the precubical set $K$ such that $K_0=K_1$ is a singleton and $K_n=\varnothing$ for $n\geq 2$: $K$ is a loop. Then $K$ is spatial but not proper (see \cite[Example~(1.5)]{MR3722069}). It means that the left-hand inclusion is strict. 

Consider the precubical set $K = \square[3] \sqcup_{\de\square[3]}\square[3]$. Any tame natural $d$-path from $0_3$ to $1_3$ lying in the common boundary $|\de\square[3]|_{geom}$ that does not contain other vertices than $0_3$ and $1_3$ is represented by two distinct elements of $\P_{0_3,1_3}|K|_{nat}$. Thus, $\square[3] \sqcup_{\de\square[3]}\square[3]$ is not spatial. It means that the right-hand inclusion is strict.
\epf

It is necessary to recall a basic fact about $\Delta$-inclusions before proceeding to the proof of Theorem~\ref{main2}. 

\bp \label{DeltaHomeo2} \cite[Proposition~2.2 and Corollary~2.3]{Moore2}
A continuous bijection $f:U\to V$ of $\top$ is a homeomorphism if and only if it is a $\Delta$-inclusion, i.e. a set map $[0,1]\to A$ is continuous if and only if the composite set map $[0,1]\to A\to B$ is continuous.
\ep

Theorem~\ref{main2} states intuitively that the $\Delta$-generated spaces $\P_{\alpha,\beta}|K|_{nat}$ for $(\alpha,\beta)$ running over $K_0\p K_0$ do not contain too many open subsets when $K$ is a spatial precubical set. 

\bth \label{main2}
Let $K$ be a spatial precubical set. The natural map of flows $|K|_{nat}\to |K|_{tc}$ is an isomorphism. In particular, for all $(\alpha,\beta)\in K_0\p K_0$, the continuous bijection $\P_{\alpha,\beta}|K|_{nat}\to \P_{\alpha,\beta}|K|_{tc}$ is a homeomorphism. 
\eth

\bpf The map of flows $|K|_{nat}\to |K|_{tc}$ induces the identity on $K_0$. For all $(\alpha,\beta)\in K_0\p K_0$, the space $\P_{\alpha,\beta}|K|_{nat}$ is Hausdorff, the space $\P_{\alpha,\beta}|K|_{tc}$ being Hausdorff. Consider a set map $f:[0,1] \longrightarrow\P_{\alpha,\beta}|K|_{nat}$ such that the composite set map $[0,1] \longrightarrow \P_{\alpha,\beta}|K|_{nat}\longrightarrow \P_{\alpha,\beta}|K|_{tc}$ is continuous. Since $[0,1]$ is path-connected, there exists a commutative diagram of spaces of the form 
\[
\xymatrix@C=2em@R=2em
{
	[0,1] \ar@{=}[d]\fr{}& \P_{\alpha,\beta}|K|_{nat}\to \P_{\alpha,\beta}|K|_{tc}\fr{}& \displaystyle\coprod\limits_{n\geq 1} \ttop([0,n],|K|_{geom})\ar@{=}[d]\\
	[0,1] \fr{}& \ttop([0,n_0],|K|_{geom}) \fr{}& \displaystyle\coprod\limits_{n\geq 1} \ttop([0,n],|K|_{geom})
}
\]
for some integer $n_0\geq 1$. By adjunction, we obtain a continuous map $[0,1] \p [0,n_0] \longrightarrow |K|_{geom}$. Since the topological space $|K|_{geom}$ is a CW-complex, the image of the compact $[0,1] \p [0,n_0]$ is a closed compact subset of $|K|_{geom}$ which, by \cite[Proposition~A.1]{MR1867354}, intersects N interiors of cubes and vertices with $N >0$ finite. We want to prove that the set map $f:[0,1] \to \P_{\alpha,\beta}|K|_{nat}$ is continuous. Since the $\Delta$-generated spaces are sequential, it suffices to prove the sequential continuity of $f:[0,1] \to \P_{\alpha,\beta}|K|_{nat}$. Let $(t_k)_{k\geq 0}$ be a sequence of $[0,1]$ which converges to $t_\infty\in [0,1]$. For all $t\in [0,1]$, $\carrier(f(t))$ is of the form $(c_1,\dots,c_p)$ with $\dim(c_i)+\dots +\dim(c_p)=n_0$. We deduce that the set $\{\carrier(\gamma(t))\mid t\in [0,1]\}$ has at most $(N+1)^{n_0}$ elements, i.e. that it is finite. Thus the sequence of carriers $(\carrier(f(t_k))_{k\geq 0}$ has a constant subsequence. Suppose that the sequence $(\carrier(f(t_k))_{k\geq 0}$ is constant and equal to $(c_1,\dots,c_p)$. Then the sequence of paths $(f(t_k))_{k\geq 0}$ belongs to the image of the continuous map $N_{\dim(c_1)}\p \dots \p N_{\dim(c_p)} \to \P_{\alpha,\beta}||K||\iso \P_{\alpha,\beta}|K|_{nat}$. The product $N_{\dim(c_1)}\p \dots \p N_{\dim(c_p)}$ is a finite product in $\top$ of compact metrizable spaces by Proposition~\ref{Delta}. By \cite[Lemma~6.9]{Moore2}, this product coincides with the product taken in the category of general topological spaces. It means that $N_{\dim(c_1)}\p \dots \p N_{\dim(c_p)}$ is compact metrizable, and hence sequentially compact. Consequently, the image of $N_{\dim(c_1)}\p \dots \p N_{\dim(c_p)} \to \P_{\alpha,\beta}|K|_{nat}$ is sequentially compact, and also closed in $\P_{\alpha,\beta}|K|_{nat}$, the latter being Hausdorff. We deduce that the sequence $(f(t_k))_{k\geq 0}$ has a limit point which is necessarily $f(t_\infty)$ by continuity of the composite map $[0,1] \longrightarrow \P_{\alpha,\beta}|K|_{nat}\longrightarrow \P_{\alpha,\beta}|K|_{tc}$. In fact, we have proved that every subsequence of $(f(t_k))_{k\geq 0}$ has a subsequence which has a limit point which is necessarily $f(t_\infty)$. Suppose that the sequence $(f(t_k))_{k\geq 0}$ does not converge to $f(t_\infty)$. Then there exists an open neighborhood $V$ of $f(t_\infty)$ in $\P_{\alpha,\beta}|K|_{nat}$ such that for some $M\geq 0$, and for all $k\geq M$, $f(t_k)\in V^c$, the complement of $V$, the latter being closed in $\P_{\alpha,\beta}|K|_{nat}$. Thus, $(f(t_k))_{k\geq M}$ cannot have a limit point: contradiction. We deduce that $f:[0,1] \longrightarrow\P_{\alpha,\beta}|K|_{nat}$ is sequentially continuous, hence continuous. It means that the continuous bijection $\P_{\alpha,\beta}|K|_{nat} \to \P_{\alpha,\beta}|K|_{tc}$ is a $\Delta$-inclusion. Therefore, the latter is a homeomorphism by Proposition~\ref{DeltaHomeo2}. The proof is complete.
\epf

Note that the functor $|\!-\!|_{tc}:\square^{op}\set\to \dtop$ is not colimit-preserving. Otherwise, there would be an isomorphism $|K|_{nat}\to |K|_{tc}$ for all precubical sets $K$, which would contradict Proposition~\ref{referee}.

As already noticed at the very end of Section~\ref{natural_sec}, the spaces of $d$-paths of precubical sets are equipped in \cite{MR4070250} with the compact-open topology instead of some kind of kelleyfication of the compact-open topology. The latter is the correct internal hom, both for $k$-spaces and $\Delta$-generated spaces, except in very specific situations like Proposition~\ref{Delta} and Proposition~\ref{dDelta}. Since the $\Delta$-kelleyfication functor takes (weak resp.) homotopy equivalences to (weak resp.) homotopy equivalences, this point is not an issue.

\bth \label{final}
	Let $K$ be a precubical set. The natural map of flows $|K|_{nat}\to |K|_{tc}$ is a weak equivalence of the h-model structure of flows. In particular, for all $(\alpha,\beta)\in K_0\p K_0$, the continuous map $\P_{\alpha,\beta}|K|_{nat}\to \P_{\alpha,\beta}|K|_{tc}$ is a homotopy equivalence.
\eth

\bpf
The map of flows $|K|_{nat}\to |K|_{tc}$ induces the identity on $K_0$. Let $(\alpha,\beta)\in K_0\p K_0$. Let $\mathbf{c}=\square[\underline{n}]\to K$ be an object of $\Ch_{\alpha,\beta}(K,n)$ with $\ell(\underline{n})=p$ and $\underline{n}=(n_1,\dots,n_p)$. Since there is an isomorphism of flows $||\square[\underline{n}]|| \iso |\square[\underline{n}]|_{nat}$ by Theorem~\ref{iso-cube}, there is the homeomorphism $N_{n_1}\p \dots \p N_{n_p} \iso \P_{0_{n_1},1_{n_p}} |\square[\underline{n}]|_{nat}$, the identity map $\square[\underline{n}]\to \square[\underline{n}]$ being the final object of the small category $\Ch_{0_{n_1},1_{n_p}}(\square[\underline{n}],n)$ (the proof is similar to the one of Proposition~\ref{lem_prep}). The precubical set $\square[\underline{n}]$ is spatial by Proposition~\ref{referee}, being proper. Thus the topological space \[\D_{\alpha,\beta}(K,n)(\mathbf{c}) = N_{n_1}\p \dots \p N_{n_p} \iso \P_{0_{n_1},1_{n_p}} |\square[\underline{n}]|_{nat}\] is homeomorphic to the space of tame natural $d$-paths in $|\square[\underline{n}]|_{geom}$ from the initial state $0_{n_1}$ of $\square[n_1]$ to the final state $1_{n_p}$ of $\square[n_p]$ equipped with the $\Delta$-kelleyfication of the relative topology induced by the compact-open topology by Theorem~\ref{main2}. It implies that the map $\P_{\alpha,\beta}||K||\to \P_{\alpha,\beta}|K|_{tc}$ is the homotopy equivalence $\coprod_{n\geq 1} F_n^K$ of \cite[Proposition~9.7]{MR4070250}. The proof is complete thanks to the isomorphism of flows $||K|| \iso |K|_{nat}$.
\epf

There exists a continuous bijection between Hausdorff $\Delta$-generated spaces which is a homotopy equivalence and not a homeomorphism: consider $X=\mathbf{S}^1$, the discretization $X^\delta$ and the map between unreduced cones $C(X^\delta)\to C(X)$ \cite{376474}. Consequently, it is not possible to deduce Theorem~\ref{main2} from Theorem~\ref{final} and Proposition~\ref{prespatial} without any additional assumption.

The word ``weak homotopy equivalence'' can be replaced by ``homotopy equivalence'' in the statements of \cite[Theorem~7.5 and Theorem~7.6]{MR4070250} because all maps of \cite[Equation~7.5]{MR4070250} are homotopy equivalences. Indeed, it is proved in \cite[Proposition~10.3]{MR4070250} that some specific map is an h-cofibration. Therefore the diagram of \cite[Proposition~10.4]{MR4070250} is Reedy h-cofibrant and the map $Q^K_n$ is a homotopy equivalence. 

\begin{cor} \label{app_general}
Let $K$ be a precubical set. Let $\alpha,\beta$ be two vertices of $K$. Then the space of execution paths $\P_{\alpha,\beta}|K|_q$ is homotopy equivalent to the space of nonconstant $d$-paths from $\alpha$ to $\beta$ in the geometric realization of $K$ equipped with the $\Delta$-kelleyfication of the compact-open topology. 
\end{cor}

\bpf
It is a consequence of Theorem~\ref{natural-equal}, Theorem~\ref{final} and \cite[Theorem~7.6]{MR4070250}. 
\epf

Since the natural realization functor is an m-realization functor, Theorem~\ref{final} provides a model category explanation of the known fact that the space of nonconstant $d$-paths in the geometric realization of a precubical set between two vertices of the precubical set has the homotopy type of a CW-complex.

\appendix

\section{Characterization of spatial precubical sets}
\label{spatial_app}

The purpose of this appendix is to characterize spatial precubical sets without using any realization functor from precubical sets to flows. This result is unnecessary for the understanding of the core of the paper. It is the reason why it is expounded here. 

\begin{nota}
	Let $n\geq 1$. Let $\mathcal{B}_n$ be the set of precubical sets $A$ such that $A\subset \de\square[n]$ and such that $|A|_{geom} \subset [0,1]^n$ contains a natural $d$-path from $0_n$ to $1_n$ in $[0,1]^n$ which does not intersect $\{0,1\}^n\backslash\{0_n,1_n\}$. In particular, it means that $0_n,1_n$ are two vertices of $A$. One has $\mathcal{B}_1 = \mathcal{B}_2=\varnothing$ and for all $n\geq 3$, one has $\de\square[n]\in \mathcal{B}_n$.
\end{nota}

Since every $d$-path of $[0,1]^n$ has a naturalization \cite[Definition~2.14]{MR2521708}, the adjective ``natural'' is superfluous in the definition of $\mathcal{B}_n$.

\bth \label{carac_spatial}
The class of spatial precubical sets is a small orthogonality class. More precisely, it is the class of precubical sets which are orthogonal with respect to the set of maps of precubical sets \[\bigg\{\square[n]\sqcup_A \square[n]\longrightarrow \square[n]\mid n\geq 3 \hbox{ and }A\in \mathcal{B}_n\bigg\}.\]
\eth

Every map $\square[n] \sqcup_{A} \square[n] \to \square[n]$ for $n\geq 3$ and $A\in \mathcal{B}_n$ being an epimorphism of precubical sets, injective is equivalent to orthogonal. Recall that the injectivity (orthogonality resp.) condition means that any map $\square[n]\sqcup_A \square[n]\to K$ factors (uniquely resp.) as a composite map $\square[n]\sqcup_A \square[n]\to \square[n] \to K$ \cite[Definition~4.1 and Definition~1.32]{TheBook}.

\bpf
Let $K$ be a precubical set which is not injective with respect to $\square[n] \sqcup_{A} \square[n] \to \square[n]$ for some $n\geq 3$ and some $A\in \mathcal{B}_n$. It means that there exists a map of precubical sets $f:\square[n] \sqcup_{A} \square[n] \to K$ which does not factor as a composite $\square[n] \sqcup_{A} \square[n] \to \square[n]\to K$. By definition of $\mathcal{B}_n$, there exists a tame natural $d$-path $\gamma$ from $0_n$ to $1_n$ such that $\gamma([0,n])\subset |A|_{geom} \subset [0,1]^n$ and such that $\gamma([0,n])\cap \{0,1\}^n=\{0_n,1_n\}$. Let
\[
\xymatrix@C=6em
{
	c_1\sqcup c_2:\square[n] \sqcup \square[n] \fr{\id_{\square[n]}\sqcup \id_{\square[n]}} & \square[n] \sqcup_{A} \square[n] \fr{f} & K.
}
\]
One obtains a composite map of spaces \[\P_{0_{n},1_{n}}|\square[n] \sqcup_{A} \square[n]|_{nat} \longrightarrow \P_{f(0_n),f(1_n)} |K|_{nat} \longrightarrow \ttop([0,n],|K|_{geom})\] such that $\xi_k\in \P_{f(0_{n}),f(1_{n})} ||K||$ represented by $\gamma\in \D_{f(0_{n}),f(1_{n})}(K,n)(c_k) = N_n$ for $k=1,2$ is taken to the natural $d$-path $[c_1;\gamma]=[c_2;\gamma]$ in $|K|_{geom}$. We have $\xi_1\neq \xi_2$ in the colimit $\P_{f(0_{n}),f(1_{n})} ||K||$ (it is the same argument as in the proof of Proposition~\ref{referee}). It means that $K$ is not spatial. 

Conversely, let $K$ be a precubical set which is injective with respect to $\square[n] \sqcup_{A} \square[n] \to \square[n]$ for all $n\geq 3$ and all $A\in \mathcal{B}_n$. Let $(\alpha,\beta)\in K_0\p K_0$. Let $\xi_1$ and $\xi_2$ be two execution paths of $\P_{\alpha,\beta}|K|_{nat}$ which are taken to the same tame natural $d$-path $\gamma = [c_1;\gamma_1]*\dots *[c_p;\gamma_p]$ in $|K|_{geom}$ with $p\geq 1$ and $0=t_0<  t_1 < \dots < t_p=\ell$ such that for all $1\leq i\leq p$ and $t\in [t_{i-1},t_i]$, $\gamma(t)=[c_i;\gamma_i(t)]$ with $\dim(c_i)\geq 1$, $\gamma([0,1])\cap K_0=\{\gamma(t_i)\mid 0\leq i \leq p\}$ and $\gamma(]t_{i-1},t_i[) \cap K_0 = \varnothing$ for $1\leq i \leq p$. This implies that $\gamma_i(]t_{i-1},t_i[) \cap \{0,1\}^{\dim(c_i)} = \varnothing$ for $1\leq i \leq p$. Let $m=\sum_i \dim(c_i)$ and  $\underline{m}=(\dim(c_1),\dots,\dim(c_p))$. For $k\in \{1,2\}$, there exists a map of precubical sets $a_k: \square[\underline{m}] \to K$ such that \[(\gamma_1,\dots,\gamma_p)\in \D_{\alpha,\beta}(K,m)(a_k)=N_{\dim(c_1)}\p \dots \p N_{\dim(c_p)}\] is identified to $\xi_k$ in the colimit $\P_{\alpha,\beta}||K||$. Let $a_k=a_k^1*\dots *a_k^p$ for $k=1,2$ with $a_k^i:\square[\dim(c_i)]\to K$. 

Choose $i \in \{1,\dots,p\}$. Let $G$ be the set of $p$-cubes $c$ of $\square[\dim(c_i)]$ with $p\geq 1$ such that $\gamma_i(]t_{i-1},t_i[)$ intersects $|c|_{geom}(]0,1[^{\dim(c)})$, say in $\gamma_i(t_c)$. By hypothesis, there are the equalities \[[a_1^i;\gamma_i(t_c)]=[a_2^i;\gamma_i(t_c)]=[c_i;\gamma_i(t_c)]\] in $|K|_{geom}$. It implies that  \[[c_i;\gamma_i(t_c)]\in |a_1^i(c)|_{geom}(]0,1[^{\dim(c)}) \cap |a_2^i(c)|_{geom}(]0,1[^{\dim(c)}).\] However, there is a bijection of \textit{sets} 
\begin{equation} \label{partition} \tag{P}
	|K|_{geom} \iso K_0 \sqcup \coprod_{p\geq 1} \coprod_{e\in K_p} |e|_{geom}(]0,1[^p).
\end{equation}
It implies that for all $c\in G$, there is the equality $a_1^i(c)=a_2^i(c)$. Since $a_1^i,a_2^i:\square[\dim(c_i)]\to K$ are maps of precubical sets, we deduce that $a_1^i(c)=a_2^i(c)$ for all $c\in G$ and all their iterated faces, i.e. on the cubes $c$ of the precubical set $\widehat{G}\subset \square[\dim(c_i)]$ generated by $G$. Let $x\in \gamma_i(]t_{i-1},t_i[)$. By \eqref{partition} applied to the precubical set $\square[\dim(c_i)]$ and since $\gamma(]t_{i-1},t_i[) \cap \{0,1\}^{\dim(c_i)} = \varnothing$, there exist $p\geq 1$ and a $p$-cube $c$ of $\square[\dim(c_i)]$ such that $x\in |c|_{geom}(]0,1[^{\dim(c)})$. Such a $p$-cube $c$ thus belongs to $G$. It implies that $\gamma_i([t_{i-1},t_i])\subset |\widehat{G}|_{geom}$. 

There are two mutually exclusive cases: $\id_{[\dim(c_i)]}\in G$ and $\id_{[\dim(c_i)]}\notin G$. In the first case, there is the equality $\widehat{G}=\square[\dim(c_i)]$. We obtain $a_1^i=a_2^i$. In the second case, there is the inclusion $\gamma_i([t_{i-1},t_i]) \subset |\widehat{G}|_{geom} \subset |\de\square[\dim(c_i)]|_{geom}$. Since $\gamma(]t_{i-1},t_i[) \cap \{0,1\}^{\dim(c_i)} = \varnothing$, we deduce that $\widehat{G}\in \mathcal{B}_{\dim(c_i)}$ by definition of the latter set. It means that the map of precubical sets $a_1^i \sqcup a_2^i : \square[\dim(c_i)] \sqcup \square[\dim(c_i)] \to K$ factors as a composite \[\square[\dim(c_i)] \sqcup \square[\dim(c_i)] \longrightarrow \square[\dim(c_i)] \sqcup_{\widehat{G}} \square[\dim(c_i)]\longrightarrow K,\] and therefore as a composite \[\square[\dim(c_i)] \sqcup \square[\dim(c_i)] \longrightarrow \square[\dim(c_i)] \sqcup_{\widehat{G}} \square[\dim(c_i)] \longrightarrow \square[\dim(c_i)] \to K.\] It implies that $a_1^i=a_2^i$ for all $i \in \{1,\dots,p\}$.

We deduce that $a_1=a_2:\square[\underline{m}] \to K$ are the same map of precubical sets. Thus $\xi_1=\xi_2$ in the colimit $\P_{\alpha,\beta}||K||$. It means that $K$ is spatial.  
\epf

\begin{cor}
	Let $K$ be a precubical set of dimension $2$ (i.e. $K_n=\varnothing$ for all $n\geq 3$). Then $K$ is spatial. 
\end{cor}

\bpf
There are no maps from $\square[n]\sqcup_A \square[n]$ to $K$ for all $n\geq 3$ and all $A\in \mathcal{B}_n$.
\epf

\let\oldaddcontentsline\addcontentsline% Store \addcontentsline
\renewcommand{\addcontentsline}[3]{}% Make \addcontentsline a no-op

%\bibliographystyle{../plainurlwithoutprefixDOI} 
%%\bibliographystyle{plain} 
%
%{%\footnotesize
%	\bibliography{../Bibliotheque}
%}

\let\addcontentsline\oldaddcontentsline% Restore \addcontentsline
\end{document}